\documentclass[3p]{elsarticle}

\usepackage[dvipsnames]{xcolor}
\usepackage[unicode=true,bookmarks=true,bookmarksnumbered=false,bookmarksopen=false, breaklinks=false,hidelinks,colorlinks=true,hypertexnames=false]{hyperref}
\hypersetup{pdfauthor={Mustafa A. Mohamad, Themistoklis P. Sapsis},pdftitle={Probabilistic response and rare events in   Mathieu's equation  under correlated parametric excitation}}
\hypersetup{colorlinks=true, linkcolor=Magenta,citecolor=Magenta, urlcolor=Magenta}

\usepackage{caption}
\usepackage{subcaption}

\journal{Ocean Engineering}
\usepackage{lmodern}
\usepackage[kerning=true,tracking=true,spacing=true,stretch=10,shrink=10]{microtype}

\usepackage{amssymb}

\usepackage{color}

\usepackage{mathtools}
\newcommand{\overbar}[1]{\mkern 1mu\overline{\mkern-1mu#1\mkern-1mu}\mkern 1mu}
\DeclarePairedDelimiter\abs{\lvert}{\rvert}
\DeclarePairedDelimiter\norm{\lVert}{\rVert}

\DeclareMathOperator{\prob}{\mathbb P}
\DeclareMathOperator{\pdf}{\rho}
\newcommand{\blank}{\,\cdot\,}

\usepackage{cleveref}
\begin{document}

\title{Probabilistic response and rare events in Mathieu's equation under correlated parametric excitation}


\author{Mustafa A. Mohamad}
\ead{mmohamad@mit.edu}

\author{Themistoklis P. Sapsis\corref{mycorrespondingauthor}}
\cortext[mycorrespondingauthor]{Corresponding author}
\ead{sapsis@mit.edu}

\address{Department of Mechanical Engineering, Massachusetts Institute of Technology,\\ 77 Massachusetts Ave, Cambridge, MA 02139, USA}

\begin{abstract} 

We derive an analytical approximation to the probability distribution function (pdf) for the response of   Mathieu's equation under parametric excitation by a random process  with a spectrum peaked at the main resonant frequency, motivated by the problem of large amplitude ship roll resonance in random seas.  The inclusion of random stochastic excitation    renders the otherwise straightforward response to a system undergoing \emph{intermittent resonances}:  randomly occurring large amplitude bursts. Intermittent resonance occurs precisely when the random parametric excitation pushes the system into the instability region, causing an extreme magnitude   response. As a result, the statistics are characterized by heavy-tails.  We apply a recently developed mathematical technique, the  probabilistic decomposition-synthesis method,   to derive an analytical approximation to the  non-Gaussian pdf of the response.  We illustrate the validity of this analytic approximation  through comparisons with Monte-Carlo simulations that  demonstrate  our result accurately captures the strong non-Gaussianity of the response. 
\end{abstract}

\maketitle


\section{Introduction}
\begin{figure}[htb]
    \centering
    \includegraphics[width=0.5\textwidth]{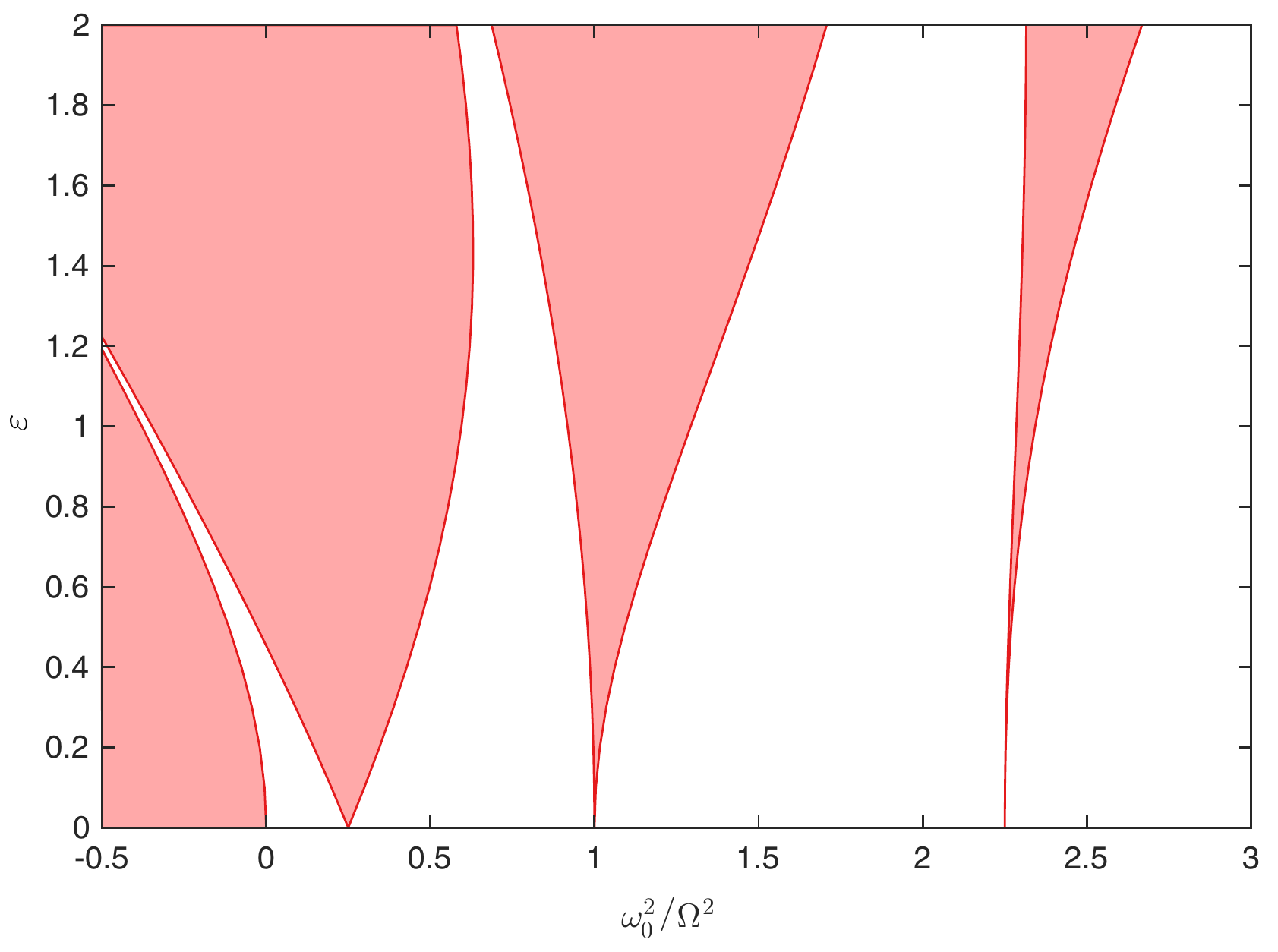}
    \caption{Ince-Strutt diagram showing the classical resonance tongues (shaded) for particular combinations of $\varepsilon$ and $\omega_0^2/\Omega^2$, for    Mathieu's equation~\labelcref{eq:simple_mat}. Transition frequencies at $\omega_0^2/\Omega^2= (n/2)^2$. Stability boundaries are  computed  using Hill's infinite determinant approach~\cite{nayfeh_mook}.}
    \label{fig:strutt}
\end{figure}

Parametrically forced systems arise in many engineering systems  and natural phenomena. 
For such systems, parametric (subharmonic) resonance can produce a large   response even when the parametric excitation amplitude is small.  To investigate these ideas, parametric resonance has been extensively  studied using  the classic  Mathieu's equation and the more general Hill's equation,
\begin{equation}\label{eq:simple_mat}
\ddot x(t) + (\omega_0^2  +\varepsilon\Omega^{2} \sin(\Omega t)) x(t) = 0,
\end{equation}
since it represents the typical response of a system when excited by a time periodic force.  When the system's natural frequency and the excitation frequency are near  the ratio $\omega_0 \,{:}\, \Omega = n \,{:}\, 2,$ for positive integers $n$, the parametric  resonance phenomena is important   and we have regions of instabilities (\cref{fig:strutt}). The most prominent resonance region occurs  for $n=1$ with the  ratio $1\,{:}\,2$. \Cref{eq:simple_mat} is a simple, deterministic  system  that is often used to    model the typical effect of parametric excitations that are  time periodic or nearly time periodic. However, this is only a special case of   more typical aperiodic time variations; indeed, in many physical systems, the induced excitations are inherently random (e.g. ship motions in random water waves~\cite{Chai2015, Lin1995, Kougioumtzoglou2014, Kreuzer}, modes in turbulence~\cite{Majda2014a, tong15, qi15}, and beam buckling due to random axial and lateral forcing~\cite{Abou-Rayan1993, Lin_Cai95}). This randomness can significantly alter the system   response,  and in certain parametric regimes   leads to   transient instabilities and intermittent bursts of extreme magnitude. These events are directly connected with the finite correlation  time of the excitation processes, and therefore cannot be quantified using deterministic analysis or approximations of the excitation processes by white-noise. 

\subsection{Motivation}\label{sec:motivation}
One important parametrically excited system in ocean engineering, and our motivating example, is ship rolling in the presence of random seas~\cite{belenky07, Chai2015, Lin1995, Kougioumtzoglou2014, Kreuzer,arnold_l,mit.00201659320120101}. The first modern theoretical description of parametric rolling was given in~\cite{Paulling59}  and    the first experimental observation of the phenomenon in~\cite{paulling1975ship}. For  ship rolling, the parametric resonance phenomena is a considerable threat to the safety of a vessel.  Indeed, ever since the investigation into the  post-Panamax C11 class cargo ship accident on October 20th, 1998, which was confirmed to have undergone severe roll motions in head seas during a storm, causing extensive loss and damage to the vessel (see~\cite{000180174200001n.d.} for details), interest in the problem has been renewed. This has lead  to developments of important guidelines and   criteria    assessing  the risk and susceptibility of vessels to roll motions~\cite{ABS2004}. As such, the problem of  parametric rolling has been an important factor in the present debate on the second generation intact stability criteria (see~\cite{Francescutto15} for the current status of the IMO criteria on rolling).

It is  now well understand that  large amplitude roll motions can occur through parametric  resonance,  even when there are no direct wave-induced roll  moments. This is most prominent in head or following seas, where roll-restoring characteristics can vary   significantly in time  compared to still water conditions. When a wave crest is amid-ship, stability is reduced as the   bow and stern are likely to have emerged, which reduce   roll-restoring moments. The effects are most pronounced for wave lengths comparable  to the ship length and increase for steeper waves. 

The roll motion   of a ship in following or head seas, for small pitch and heave motions, can be modeled by, 
\begin{equation}
I \ddot \phi(t) + B(\dot \phi(t)) + \Delta GZ(\phi(t)) = F(t),
\end{equation}
where $\phi(t)$ is the roll angle,  $I$ is the inertia term including added mass,  $B(\dot \phi(t))$ is the damping moment, $F(t)$ is a small wave excitation term, and $\Delta GZ(\phi(t))$ is the restoring moment term. The time-dependent restoring moments in irregular seas may, sometimes, be approximated by~(see e.g. \cite{belenky07,Kreuzer})
\begin{equation}
\Delta GZ(\phi(t)) = (\alpha + \xi(t))\phi(t) - \beta \phi^3(t),
\end{equation}
where $\xi(t)$ is the random wave excitation,  with a narrowbanded spectrum $S_\xi(\omega)$ around a frequency  $\Omega$; thus we expect the parametric resonance phenomena to be important when the encounter frequency $\Omega$ is roughly twice the natural  natural frequency of rolling.  In addition, $B(\dot\phi(t))$ is    a nonlinear function of the roll   velocity. Since our goal  is to model the leading order probabilistic dynamics that occur due to the interaction of the random  time-dependent restoring  term  and the vessel's roll angle, we  neglect the assumed small nonlinear terms in the restoring   and damping moments:
\begin{equation}\label{eq:ship_roll_prototype}
  \phi(t) + f \dot\phi(t)  +(\alpha + \xi(t))\phi(t)  = F(t);
\end{equation}
this choice is motivated by the fact that stochastic roll resonance, in this context, is a consequence of the multiplicative excitation term,  and not nonlinearities in the restoring or damping terms.  However, it is worthwhile to remark on  the theoretical nature of~\cref{eq:ship_roll_prototype}; system nonlinearities are important  in   accurate models of ship rolling and play a role in modifying the instability zones~\cite{Neves20071618}. The inclusion of the cubic term in the  restoring moment would, for one,    impact the underlying shape at  the very tail ends  of the response  probability density function  by suppressing their magnitude. Furthermore,  the excitation in roll motion is  coupled   with pitch and heave motions, and this requires at least three degrees of freedom. Despite these remarks, the model in~\cref{eq:ship_roll_prototype} serves as a useful prototype system for analytical work  investigating  the complex heavy-tailed   probabilistic properties of  roll resonance in random seas; and is motivated by the desire in design to have simple and accurate  analytic predictors of dynamic stability  that  account for   extreme conditions~\cite{Sprou2000}.

This summarizes the basic equation that govern  the leading-order  ship roll dynamics in the presence of irregular seas. Finally, we note the following: in certain parametric regimes,  solutions of~\cref{eq:ship_roll_prototype}  exhibit random periods of large amplitude roll motions (intermittency). It is this parametric regime   we are interested in  investigating here. 

\subsection{Stochastic generalization of  Mathieu's equation}

While   Mathieu's equation and its variants have been extensively analyzed~\cite{dynamicsofparametricexcit,nayfeh_mook}, generalizations incorporating stochasticity  are less well understood, but have been studied in various contexts~\cite{Soong_Grigoriou93,Lin_Cai95,Poulin1885,doi:10.1137/070689322,stratonovich1967topics}. A very important  dynamical transition that   occurs in the presence of random parametric forcing is \emph{intermittent resonance}; that is, randomly occurring short-lived periods when the system experience parametric resonance. Deterministic models based on some time-averaged property of the noise   would fail to capture this   important dynamical transition, since this is an essentially   transient phenomena.   

To explore the effects of  noise in parametric forcing, motivated by~\cref{sec:motivation},  we consider the following generalization of~\cref{eq:simple_mat}:
\begin{equation}\label{eq:mat_simple_stoch}
\ddot x(t) + 2 \zeta \omega_0 \dot x(t)+  \omega_0^2(1  + \kappa(t)) x(t) =   F(t),
\end{equation}
where $\kappa(t)$ is a  narrowbanded  random process,  $F(t)$ is a broadbanded forcing term of low intensity, e.g. white noise, and  $\zeta$ is the damping
coefficient. For example, $\kappa(t)$ can be thought of as $\kappa (t) =  \alpha(t) \sin(\Omega t)$, with  frequency ratio near one of the resonance regions and $\alpha(t)$ a  random process, such that the dominant energy in the spectrum $S_{\kappa}(\omega)$ is concentrated near $\Omega$; in other words, the noise has a dominant frequency component  at $\Omega$.   We explain later in detail the motivation behind this selection, for now we note that this spectrum follows quite naturally when  modeling    physical processes, in particular, random seas.

\subsection{Perspective}
Our goal here is to derive an analytical approximation to the probability density function (pdf) for the generalized Mathieu's equation in~\cref{eq:mat_simple_stoch} when forced parametrically by a correlated random process near the principal resonance region; in particular, we are interested in   the \emph{stationary probability distribution} of the response in the regime undergoing intermittency. As mentioned, this is   challenging since the system exhibits transient  resonance, which can be identified by large amplitude spikes in the time-series of the response; as a result of intermittency, the resulting pdf is non-Gaussian with heavy-tailed characteristics.  Recently, there have been efforts to quantify the heavy-tailed statistical structures for systems undergoing intermittent instabilities~\cite{mohamad2015,mohamad2016, tong15,qi15}. Here,   we apply a recently developed technique  designed to  approximate the pdf of systems exhibiting intermittent instabilities:   the probabilistic decomposition-synthesis method~\cite{mohamad2016,mohamad2015}. Since the system we investigate is low-dimensional, with linear damping and restoring terms,  we  apply the formulation in~\cite{mohamad2015} to approximate the probability  measure of the response. This approach provides us with analytical results, \emph{taking into account the   correlated nature of the multiplicative excitation process $\kappa(t)$}. 

The benefits of this approach are numerous. For  systems, such as in~\cref{eq:ship_roll_prototype}, we can derive analytical results for the  non-Gaussian response pdf,   for both the main probability mass and the heavy-tailed structure.   For systems where nonlinearities are important, we can adapt the method and  apply the  probabilistic decomposition-synthesis method in a computational setting, as described in~\cite{mohamad2016}, for a fast approximation of the response pdf that  accounts for system nonlinearities. 

While several   methods can be applied to solve for the stationary measure for systems excited by multiplicative noise, for the case  of transient instabilities,   as might occur  in parametric rolling~\labelcref{eq:ship_roll_prototype}, they are severely limited in practicality  due to their large computational demands and/or limitations in dealing with  the strongly transient nature of intermittent instabilities.  For example, the Monte-Carlo approach (direct sampling of long time realizations  of the system), while very attractive since it  provides the most accurate results, is  a very computationally intensive procedure, requiring many realizations  for  accurately resolved  tail statistics. Another technique is to  formulate the associated Fokker-Planck-Kolmogorov (FPK) equation for~\cref{eq:mat_simple_stoch}~\cite{soize94}. This can be performed by utilizing shape filtering to  approximate the correlated excitation process.  However, using filtered Gaussian white noise is  prone to introduce significant errors in  the tails of the response (even small numerical errors in time-series simulations of $\kappa(t)$  lead to large inaccuracies in the tail statistics of $x(t)$~\cite{majda_branicki_DCDS}). In any case, this approach   requires solving a   demanding FPK equation~\cite{Masud_bergman05}, which also has to be done to high accuracy (tail events have extremely low probabilities). 

Stochastic averaging is another widely used method~\cite{Pavliotis2008}; the typical approach here   would be to first derive a set of equations for the slowly varying variables and then to  apply  the    stochastic averaging procedure to arrive at a set of  It\={o} stochastic differential equations for the transformed coordinates. The Fokker-Plank-Kolmogorov equation can then be used to solve for the response pdf~\cite{Lin_Cai95,8702927420120101}. Clearly, this approach is not valid in parametric regimes undergoing intermittent resonance, since it averages away the time dependent nature of the multiplicative excitation, leading to Gaussian statistics, which is  decidedly not the case.

\subsection{Outline}
In~\cref{sec:prob_def} we formulate the problem and provide the   problem statement; we explain  the particular form of the excitation noise structure and its interaction with the system dynamics in the parametric regime of interest. Following this discussion, in~\cref{sec:slow}, we derive equations that govern the slow dynamics of the problem. The equations for the slowly varying variables are the starting point of our application of the probabilistic decomposition-synthesis method, which we briefly give an overview in~\cref{sec:PDSmethod}. In~\cref{sec:analysis} we apply the method to derive  the  analytical formula that approximates   the pdf of Mathieu's equation in the parametric regime of interest, and in~\cref{sec:comps} we compare the analytical formula  with numerical results from Monte-Carlo simulations.

\section{Problem formulation}\label{sec:prob_def}

We consider the following stochastic generalization of   Mathieu's equation:
\begin{equation}\label{eq:canonical_stoch_Mathieu}
    \ddot x(t) + 2 \zeta \omega_0 \dot x(t) +  \omega_0^2( 1  +  \kappa(t) )  x(t) =   F(t),
\end{equation}
where  $\omega_0$ is  the undamped natural frequency of the system, $\zeta$ is the damping coefficient,  $\kappa(t)$ is a   narrowbanded  random process around $\Omega$, and   $F(t)$ is an additive  broadbanded    random excitation term; we assume both $\kappa(t)$ and $F(t)$ are stationary Gaussian processes.

As mentioned, the deterministic form of Mathieu's equation,
\begin{equation}\label{eq:Mathieu_det}
    \ddot x(t) + 2   \zeta \omega_0 \dot x(t) + \omega_0^2( 1 + \alpha \sin \Omega t )  x(t) =  0 ,
\end{equation}
has unstable solutions  depending upon  the parametric excitation frequency $\Omega$ and amplitude $\alpha$ parameters. Near $\Omega/ 2 \omega_0 = 1/n$, for positive integers $n$, we have regions of instabilities, with the widest instability region being for $n=1$. Damping has the effect of raising the instability regions from the $\Omega/ 2 \omega_0$ axis by $2 (2\zeta)^{1/n}$. Therefore, for $\zeta \ll 1$ the instability region near $n=1$ is of greatest practical importance~\cite{Lin_Cai95,nayfeh_mook}. 

\begin{figure}[htb]
    \centering
    \includegraphics[width=0.9\textwidth]{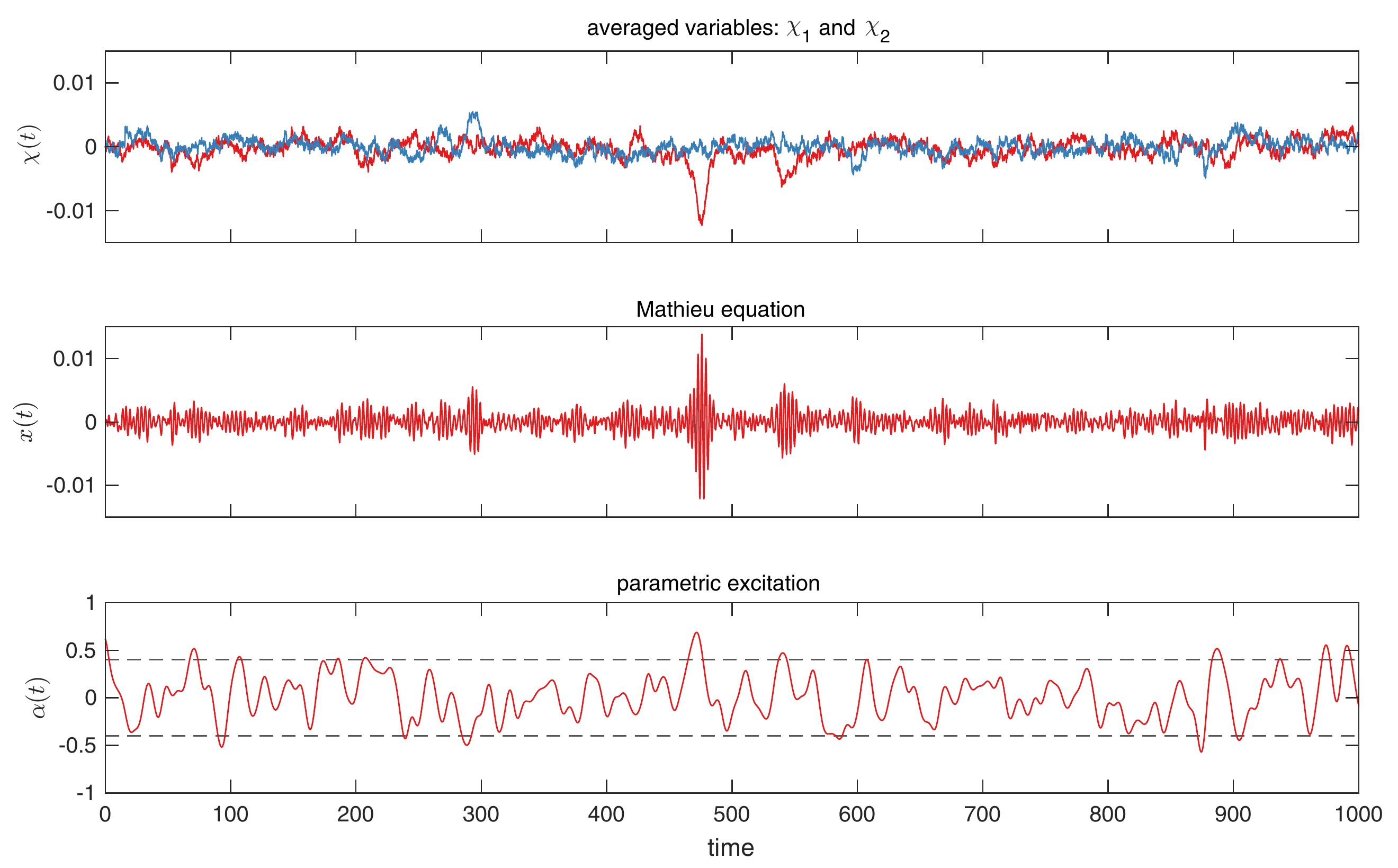}
    \caption{Sample realization of  Mathieu's equation (middle,~\cref{eq:problemstate_stoch_Mathieu}) under the prarametric excitation term (bottom) and the corresponding averaged variables $\chi_1$ and $\chi_2$ (top,~\cref{eq:averagedeqnsx1,eq:averagedeqnsx2}). The random amplitude forcing  term $\alpha(t)$  triggers intermittent resonance when it crosses above or below the instability thresholds (dashed lines).}
    \label{fig:real_mathiue}
\end{figure}

Furthermore,   since we are interested in analyzing the response pdf in the regime   where the system is undergoing intermittent resonance,  we assume $\kappa(t)$ is a narrowbanded random process around the most  important resonant  region at frequency $\Omega = 2\omega_0$. The approach can be extended if the frequency is  detuned,   but for simplicity of the presentation we consider no detuning. In this case, intermittent resonance is a possibility since the stochastic process  $\kappa(t)$  may randomly cross into the instability tongue. In other words, for the regime of interest,  on    average $\kappa(t)$ is in a state such that the variance of the response $x(t)$ is associated with the background attractor, however with low probability   $\kappa(t)$ can transition to  a critical regime, producing an intermittently  unstable  response of severe magnitude. Thus, through this regime switching, we have randomly occurring periods of large amplitude responses, which are finite-time instabilities with   positive Lyapunov exponent; a typical time-series is shown in~\cref{fig:real_mathiue}. The severity of these instabilities   depends upon the   magnitude of the damping term and the amplitude and characteristic time-scale  of the multiplicative  noise term.  

\subsection{Excitation process}
For concreteness we assume the excitation process takes a canonical form $\kappa(t) = \alpha(t) \sin(2\omega_0 t)$, where $\alpha(t)$ is a stationary Gaussian process with a non-oscillatory correlation function. Additionally,   $\overbar{\alpha(t)} = 0$ and  $R(\tau) = \overbar{\alpha(t)\alpha(t+\tau)} = \sigma_\alpha^2 e^{-\tau^2/2\ell_\alpha^2}$, where $\ell_\alpha$ is the characteristic time-scale of the process and $\sigma_\alpha^2$ its variance.


\subsection*{Problem statement}
With these considerations, the problem is to derive an approximation for the stationary probability distribution for the system:
\begin{equation}\label{eq:problemstate_stoch_Mathieu}
    \ddot x(t) + 2 \zeta \omega_0 \dot x(t) +  \omega_0^2( 1  + \alpha(t) \sin( 2 \omega_0 t) )  x(t) =   F(t),
\end{equation}
in the parametric regime undergoing intermittent resonance. The final result  is  given in~\cref{eq:mathanalyic}.

\section{Derivation of the  slow dynamics}\label{sec:slow}

We proceed by assuming a narrow band response around $\omega_0$ and averaging over this fast frequency the governing system~\labelcref{eq:problemstate_stoch_Mathieu}. By introducing the   coordinate  transformation
\begin{equation}\label{eq:transformation}
    \begin{aligned}
        x(t) &= \chi_1(t) \cos(\omega_0 t) + \chi_2(t) \sin(\omega_0 t), \\
        \dot x(t) &=  -  \omega_0 \chi_1(t) \sin(\omega_0 t) + \omega_0 \chi_2(t) \cos(\omega_0 t), 
    \end{aligned}
\end{equation}
in~\cref{eq:problemstate_stoch_Mathieu} and using the additional relation $\dot \chi_1(t) \cos(\omega_0 t) + \dot \chi_2(t) \sin(\omega_0 t) = 0$, we obtain the following pair of   equations for the slow variables  $\chi_1(t)$ and $\chi_2(t)$:
\begin{align}
    &\begin{multlined} \label{eq:diffeqTransformed1}
        \dot \chi_1   = - \biggl[  2\zeta \omega_0 \biggl(\chi_1 \sin^2(\omega_0 t) -  \frac{1}{2} \chi_2 \sin(2\omega_0 t)\biggr) -\frac{\omega_0 \alpha}{2} \chi_1 \sin^2(2 \omega_{0} t)    \\-\omega_0 \alpha \chi_2 \sin^2(\omega_0 t) \sin(2\omega _{0}t)   \biggr] - \frac{1}{ \omega_0} \sin( \omega_0 t) F(t),
    \end{multlined}\\   
    &\begin{multlined} 
        \dot \chi_2   =   \biggl[  2\zeta \omega_0 \biggl( \frac{1}{2} \chi_1 \sin(2 \omega_0 t) -    \chi_2\cos^2( \omega_0 t)\biggr)  -\frac{\omega_0 \alpha}{2} \chi_2 \sin^2(2 \omega_{0} t)     \\-\omega_0 \alpha  \chi_1 \cos^2(\omega_0 t) \sin(2\omega_{0} t)   \biggr]  + \frac{1}{ \omega_0} \cos( \omega_0 t) F(t).
    \end{multlined} \label{eq:diffeqTransformed2}
\end{align}
Averaging the deterministic terms in brackets over the fast frequency  in~\cref{eq:diffeqTransformed1,eq:diffeqTransformed2}   gives,
\begin{align}\label{eq:averagedeqnsx1}
    \dot \chi_1(t) &= -\biggl( \zeta  -\frac{ \alpha(t)}{4} \biggr)\omega_0  \chi_1(t) - \frac{1}{ \omega_0} \sin( \omega_0 t) F(t),\\
    \dot \chi_2(t) &= -\biggl( \zeta  + \frac{ \alpha(t)}{4} \biggr)\omega_0  \chi_2(t) + \frac{1}{ \omega_0} \cos( \omega_0 t) F(t).\label{eq:averagedeqnsx2} 
\end{align}
The equations above for the slow variables provide good pathwise and statistical agreement with~\cref{eq:problemstate_stoch_Mathieu}. Furthermore, these equations make it clear that an instability is expected when  $\abs{\alpha(t)} > 4\zeta$ (shown in dashed lines in~\cref{fig:real_mathiue}); this instability criterion is nothing more than the    well-known instability condition of the deterministic case (for fixed in time $\alpha$), which   to leading order is given by: $\delta = 1/4 \pm 1/2\sqrt{\alpha^2\delta^2-4\zeta^2\delta}$, where $\delta = \omega_0^2/\Omega^2$ (see left plot in~\cref{fig:superpos})

Next, we apply a stochastic averaging procedure to the additive forcing term, also known as the diffusion approximation~\cite{Lin_Cai95, Klyatskin05, sapsis_PEM11}. More specifically, if the governing dynamics  act on a sufficiently slower time scale than the memory of the additive stochastic process, then the independent increment approximation is valid. This is the case if the stochastic process $F(t)$ is broadbanded, and leads to the following  set of It\={o} stochastic differential equations for the slow variables:
\begin{align}\label{eq:ouProcAvedx1}
    \dot \chi_1(t) &= -\biggl( \zeta  -\frac{ \alpha(t) }{4} \biggr)\omega_0  \chi_1(t) +\sqrt{2 \pi K } \, \dot W_1(t),\\
    \dot \chi_2(t) &= -\biggl( \zeta  + \frac{ \alpha(t)}{4} \biggr)\omega_0  \chi_2(t) +\sqrt{2 \pi K } \, \dot W_2(t),\label{eq:ouProcAvedx2}
\end{align}
with $K =S_{F}(\omega_0)/ 2 \omega_0^2$, where $S_{F}(\omega_0)$ is the spectral density of the additive excitation $F(t)$ at frequency $\omega_0$, and $\dot W_1$ and $\dot W_2$ are independent white noise processes of unit intensity~\cite{Lin_Cai95}. The slowly varying variables, after averaging the forcing term, transform to two decoupled stochastic differential equations. \Cref{eq:ouProcAvedx1,eq:ouProcAvedx2} are a good statistical approximation  to the original system~\labelcref{eq:problemstate_stoch_Mathieu} (but provide poor pathwise agreement). 

We emphasize that the broadband hypothesis for the additive stochastic process is a convenient setup that leads to the derived white-noise formulation. However, the analysis and results that follow do not require the white-noise formulation and are directly applicable for the general case of an (non white-noise) additive stochastic forcing. Such a case could be, for example, an additive narrowbanded forcing term with spectral content distributed around $\omega_0$. 

\section{The probabilistic decomposition-synthesis method}\label{sec:PDSmethod}

The analysis in~\cref{sec:prob_def} provides  the dynamics of the slow variables and  reveals the  interaction of  the parametric excitation process $\alpha(t)$ with the slow variables. Starting from these equations we can apply the  probabilistic decomposition-synthesis method to analytically approximate the stationary measure of $\chi(t)$. To be self-contained,  we provide a very brief overview  of the probabilistic decomposition-synthesis method adapted to the current problem; further details can be found in~\cite{mohamad2015}  and a detailed description of the method in a general context in~\cite{mohamad2016}.

The main idea of the method is to decompose the system response as:
\begin{equation}
x(t)=x_b(t)+x_r(t),
\end{equation}
where $x_r$  is the solution when a rare event due to an instability occurs and $x_b$ is the stochastic response otherwise, i.e. the background dynamics. To be more specific  $x_r$ will be the response of the system when the following two conditions are satisfied: (i) $\norm{x}>\zeta$, where $\zeta$ is the extreme event threshold with respect to a chosen norm $\norm{\blank}$, and (ii) the parametric excitation  obtains values that lead to an instability, i.e. $\alpha(t)\in R_e$, where $R_e$ describes the critical region of values for $\alpha(t)$ that induce  an instability.

Together with this decomposition into rare events and the stochastic background, we also adopt the following assumptions:
\begin{enumerate}
\item The existence of intermittent events has negligible effect on the statistical characteristics of the stochastic attractor and can be ignored when analyzing the background state $x_b$;
\item  Rare events are statistically independent from each other.
\end{enumerate}
With this setup we can analyze the two states separately and     probabilistically synthesize the information obtained from this  analysis. This is completed using a total probability argument to   obtain the statistics for an arbitrary quantity of interest  $q$   by
\begin{equation}\label{eq:probLaw}
\pdf(q) =  \pdf(q \mid {\norm{x}>\zeta,\alpha\in R_{e}} ) \prob_{r} +  \pdf(q \mid x=x_b) ( 1-\prob_{r} ).
\end{equation}
In this paper, the quantity of interest is the system response $x$. The first term expresses the contribution of   rare events due to   instabilities and is the heavy-tailed portion of the distribution for $q$ and the second term expresses the contribution of the background state, which contributes  the main      probability mass in the pdf for $q$. Moreover, $\prob_{r} \equiv \prob(\norm{x}>\zeta,\alpha\in R_{e}) $ is the total probability of a rare event, which is defined as the  ratio between the time   the system spends in rare event responses over the total time. Note that the temporal duration of rare transitions  also includes a decay or relaxation phase to the background attractor, where the instability is no longer active but the system  response  still has important magnitude.

\section{The probability distribution for the response}\label{sec:analysis}

Here we  apply the various steps of the probabilistic decomposition-synthesis method to derive analytical results that approximate the heavy-tailed distribution of the response for the system in~\cref{eq:problemstate_stoch_Mathieu}. In particular, we apply the method directly on the slow variables, since the fast frequency we averaged over  is  inconsequential in the    pdf of the response.

Firstly, because $\alpha(t)$ is a zero mean process both $\chi_1(t)$ and $\chi_2(t)$  follow the same  probability distribution.  Consider the following  equation that represents the slowly varying variables: 
\begin{equation}\label{eq:ouProcess}
\dot \chi(t) = -\biggl( \zeta   -\frac{ \alpha(t)}{4} \biggr)\omega_0  \chi(t)   + \sqrt{2 \pi K }  \, \dot W(t)
\end{equation}
We write~\cref{eq:ouProcess} as
\begin{equation}\label{eq:ouProcessCanonical}
\dot \chi(t) = - \gamma(t)  \chi (t)  +  \sigma_F \, \dot W(t),
\end{equation}
where $\sigma_F = \sqrt{2 \pi K } $ and $\gamma(t) =  \zeta \omega_0 -     \omega_0    \alpha(t)/ 4$ is a Gaussian process with mean $   \zeta \omega_0$ and standard deviation $ \sigma_\alpha \omega_0  / 4 $. \Cref{eq:ouProcessCanonical} will be  the starting point  for   the application of the decomposition-synthesis method.

\subsection{Decomposition and  the instability region}

 \Cref{eq:ouProcessCanonical} makes it clear  that intermittency is due to the parametric forcing term $\gamma(t)$ switching signs from positive to negative values. This sign switching causes $\chi(t)$ 
to transition from its regular response to a domain  where the likelihood of an instability is high. This   switching in $\gamma(t)$ is the   mechanism behind instabilities in the variable $\chi(t)$. Therefore, we define the instability region as
\begin{equation}
R_e \coloneqq \{ \gamma(t) \mid \gamma(t)  < 0\}.
\end{equation}
In addition, for convenience, we  define an instability  threshold  by   $\eta \coloneqq   4\zeta /   \sigma_{\alpha}$, the ratio of the mean over the standard deviation of the process $\gamma$. 

\subsection{Conditional distribution of the background dynamics}\label{sec:back_dyn}

In the background state  $R_e^c$, by definition, we have no rare events. We can approximate the background dynamics by  replacing $\gamma(t)$ with its  average value in this regime. The conditional average of $\gamma(t)$ in $R_e^c$ is  
\begin{equation}
 \overbar \gamma |_{\gamma>0} =   \zeta \omega_0 +   \frac{\sigma_{\alpha}\omega_0}{4}  \frac{\phi(\eta)}{\Phi(\eta)},
 \end{equation}
 where $\phi(\blank)$ is the normal probability density function and $\Phi(\blank)$ is the normal cumulative distribution function; thus the background dynamics is described by the Ornstein-Uhlenbeck process:
\begin{equation}\label{eq:ouProcessCanonicalStable}
\dot \chi (t) = - \overbar \gamma |_{\gamma>0} \chi(t)   +  \sigma_F \, \dot W(t).
\end{equation}
Now the dynamics are globally stable in $R_e^c$  and   we can directly obtain the stationary distribution for~\cref{eq:ouProcessCanonicalStable}:
\begin{equation}\label{eq:oustabresponse}
\pdf_{\chi}(\chi \mid R_e^c) = \sqrt{\frac{\overbar \gamma |_{\gamma>0}}{\pi \sigma_F^2 }} \exp \biggl({- \frac{\overbar \gamma |_{\gamma>0}}{\sigma_F^2}} \chi^2\biggr),
\end{equation}
which is Gaussian distributed. To formulate our result in terms of the system variable $x$, we refer to the narrowbanded approximation made when averaging the governing system~\labelcref{eq:problemstate_stoch_Mathieu}. This  gives, approximately, $x_{}=\chi \cos\varphi,$ where $\varphi$ is a uniform random variable distributed between 0 and $2\pi$. The probability density function for $z=\cos\varphi$ is given by $\pdf_{z}(z)= 1/(\pi\sqrt{1-z^{2}}),$ $z\in
[-1,1].$ To avoid additional integrations for the computation of the pdf, we approximate the distribution  for $z$ by $\pdf_{z}(z)=\frac{1}{2}(\delta(z+1)+\delta(z-1)).$
This   gives the following approximation for the background statistics
\begin{equation}\label{eq:eq_positio_unstable_oscillator}
\pdf(x \mid R_e^c)=\pdf_{\chi}(x \mid R_e^c)= \sqrt{\frac{\overbar \gamma |_{\gamma>0}}{\pi
\sigma_F^2 }} \exp \biggl({- \frac{\overbar \gamma |_{\gamma>0}}{\sigma_F^2}}
x^2\biggr).
\end{equation}Therefore, the  conditional distribution of the  background dynamics is     Gaussian distributed.
 
\subsection{Conditional distribution of  rare events}

Here we derive the conditional distribution of the response when an instability occurs. We characterize localized instabilities by a growth phase, corresponding directly to $\gamma(t) < 0$, and a relaxation phase that brings the system back to the background state; both phases follow the same distribution~\cite{mohamad2015}. Additionally, during the occurrence of an instability we neglect   additive excitation and   damping,   and approximate the magnitude of the envelope as  $\xi\equiv|\chi|\ \simeq \xi_0 e^{\Lambda T_{\gamma < 0}}$, where $\xi_0$ is  the magnitude of the position's envelope, $|\chi_0|$ right before the instability has begun to emerge, $\Lambda$ is a random variable that represents the Lyapunov exponent, and $T_{\gamma < 0}$ is the random length of time that the   process $\gamma$ spends below the zero level.

We first determine the statistical characteristics of $\Lambda $ and $T_{\gamma < 0}$ (which we assume are independent). By substituting  the representation $\xi\equiv|\chi| \simeq |\chi_0 |e^{\Lambda T_{\gamma < 0}}$ into~\cref{eq:ouProcessCanonical} we obtain $\Lambda = -\gamma$ so that
\begin{equation}\label{eq:growth_lambda}
    \pdf_\Lambda(\lambda) = \prob(-\gamma \mid \gamma<0)  = \frac{4}{\sigma_{\alpha}\omega_0(1 - \Phi(\eta)) } \phi\biggl(4\frac{\lambda +  \zeta \omega_0}{\sigma_{\alpha}\omega_0}\biggr).
\end{equation}

The distribution of the duration of time  the  process $ \gamma(t)$ spends below an arbitrary  threshold level $\eta$ is not in general available. However, one can show that the  asymptotic expression in the limit of rare crossings, $\eta \to \infty$, is~\cite{MR0094269}
\begin{equation}\label{eq:rayleigh}
\pdf_{T_{\gamma < 0}}(t)=\frac{\pi t}{2\overbar{T}^{2}}e^{-\pi t^2/4\overbar{T}^{2}},
\end{equation}
and in our case this is a   good approximation since we assume  instabilities are rare events. In~\cref{eq:rayleigh} $\overbar T$ represents the average length of an instability, which for a Gaussian process is given by the ratio between the probability of $\gamma < 0$ and the average number of downcrossings of the zero level per unit time $\overbar{{N}^{-}_\gamma}(0)$ by $\gamma$~\cite{MR0094269}
\begin{equation}
\overbar{T}=\frac{\prob( \gamma < 0)}{\overbar{{N}^{-}_\gamma}(0)} = \frac{ 1 - \Phi(\eta)}{ \frac{1}{2\pi}\sqrt{-R_{ \alpha}''(0)}\exp(-\eta^{2}/2) },
\end{equation}
where we have used Rice's formula for the expected number of upcrossings~\cite{MR0370729} and   $R_{\alpha}(\tau)$ is the correlation function of $\alpha(t)$ and $R_\alpha''(0)$ is the second  derivative of the correlation function evaluated at $\tau=0$.

With these results we can determine the distribution of $\xi$ in $R_e$; the derived distribution is given by:
\begin{equation}\label{eq:density1}
    \pdf(\xi\mid  \xi_0, \alpha \in R_e)=\frac{1}{\xi}\int\limits_{0}^{\infty}\frac{1}{y}\pdf_{\Lambda}(y)\pdf_{T}\biggl(\frac{\log(\xi/\xi_0)}{y}\biggr)\, dy, \quad \xi>\xi_0.
\end{equation}
Substituting in~\cref{eq:growth_lambda,eq:rayleigh} gives,
\begin{equation}\label{eq:densitygrowth}
    \pdf(  \xi\mid  \xi_0, \alpha \in R_e) = \frac{2\pi \log(\xi/\xi_0)}{\sigma_{\alpha} \omega_0 \overbar T^2 (1-\Phi(\eta)) \xi} \int\limits_0^\infty \frac{1}{y^2}\phi\biggl(4\frac{y + \zeta \omega_0}{\sigma_{\alpha} \omega_0}\biggr)   \exp{\bigg(-\frac{\pi}{4 \overbar T^2 y^2 } \log(\xi/\xi_0)^2\bigg)}\,dy, \quad  \xi>\xi_0.
\end{equation}
The   variable $\xi_0$  corresponds to  the magnitude of the background state   $x_b$    before an instability occurs. Since, the background state is a narrowbanded Gaussian process, the magnitude of the envelope  $\xi_0$  can be modeled as a Rayleigh distribution (see~\cite{267314419860322}) with scale parameter equal to the standard deviation of the Gaussian process~\labelcref{eq:eq_positio_unstable_oscillator}. 

Note, however, that in~\cref{eq:probLaw} we are interested in the conditional statistics of events caused by instabilities which also have important magnitude. As described in~\cite{mohamad2015}, if the envelope's magnitude when an instability occurs is   small,  then even though we have an instability we do not necessarily  have a rare event, i.e. a response that is distinguishable from the typical background
 state response. This requires us to   consider background states   $\xi_0$, which are sufficiently large to result in a rare event $\norm{x}>\zeta$. To this end, we consider only background states such that $\xi_0>\zeta= \sigma_F / \sqrt{2 \overbar \gamma|_{\gamma>0}}$, i.e. we consider   extreme   responses with magnitude at least one standard deviation of the background statistics when an instability also occurs. We emphasize that the exact choice for $\zeta$ plays   little role on the approximation properties of the derived analytical expression.  This requirement gives the final distribution for $\xi_0$,
\begin{equation}\label{eq:modeenvelope_mat}
\pdf_{\xi_0}(\xi_0\mid \xi_0 > \zeta) =  \frac{2 \overbar\gamma|_{\gamma>0}}{\sigma_F^2} (\xi_0 - \zeta)  \exp\biggl(- \frac{  \overbar\gamma|_{\gamma>0}}{\sigma_F^2} (\xi_0-\zeta)^2\biggr), \quad \xi_0 > \zeta.
\end{equation}

Using~\cref{eq:modeenvelope_mat} in~\cref{eq:densitygrowth} we obtain
\begin{equation}
\pdf_{}(\xi \mid \xi>\zeta, \alpha \in R_e) = \int \pdf(  \xi\mid  \xi_0, \alpha \in R_e)  \pdf_{\xi_0}(\xi_0\mid\xi_0 > \zeta)\,d\xi_0 .
\end{equation}
In the last step, we transform the envelope representation  back to the real variable $x$ through the same   narrowbanded argument used in~\cref{sec:back_dyn} for the  background state distribution; this gives the final result for the distribution of the rare event regime:
\begin{align}\label{eq:ouunstabregime}
 \pdf_{}(x \mid \left\Vert x\right\Vert>\zeta,\alpha \in R_e) &\simeq   \frac{1}{2} \int \pdf(\left\vert  x\right\vert\mid  \xi_0,
\alpha \in R_e)  \pdf_{\xi_0}(\xi_0\mid\xi_0 > \zeta)\,d\xi_0 ,\\
          &\begin{multlined}=  \frac{ \sqrt{ 2 \pi} \overbar\gamma|_{\gamma>0} }{\sigma_F^2 \sigma_{\alpha}\omega_0\ \overbar  T^2 (1-\Phi(\eta))   }  \int\limits_\zeta^{\abs{x}}\!\int\limits_0^\infty \frac{\log(\abs{x}/\xi_0)}{y^2(\abs{x}/(\xi_0-\zeta))} \exp\biggl(-8\frac{(y + \zeta \omega_0)^2}{ \sigma_{\alpha}^2\omega_0^2}  \\- \frac{\pi}{4 \overbar T^2 y^2 } \log(\abs{x}/\xi_0)^2 - \frac{\overbar\gamma|_{\gamma>0}}{\sigma_F^2} (\xi_0-\zeta)^2 \bigg) \,dy\,d\xi_0,\end{multlined}\nonumber
\end{align}

\subsection{Probability of rare events}

Next we determine   the total probability of rare events, that is the ratio of time that the system response spends in rare transitions over the   total time, which we denote by  $\prob_r$. This quantity can be computed by analyzing the    duration of transitions into $R_e$, i.e. the duration of instabilities, but also including the time it takes for the response to return to the background attractor.

Consider a representative extreme event with an average growth rate $\overbar{\Lambda}_+$ and decay rate $\overbar{\Lambda}_- $. During the growth phase the dynamics are approximated by:
\begin{equation}\label{eq:osc_typicalgrowth}
        \xi_p = \xi_0 \exp({\overbar{\Lambda}_+ T_{\gamma<0}}),
\end{equation}
where $T_{\gamma<0}$ is the duration of the growth event and $\xi_p$ the peak value of the response. Similarly, over the decay phase:
\begin{equation}\label{eq:osc_typicaldecay}
        \xi_0 = \xi_p \exp({-\overbar{\Lambda}_- T_\text{decay}}).
\end{equation}
This gives that: 
\begin{equation}\label{eq:Probability_connection}
\frac{T_\text{decay}}{T_{\gamma<0}} = \frac{\overbar{\Lambda}_+}{\overbar{\Lambda}_{-}} =  \frac{-\overbar{\gamma}|_{\gamma<0}}{\overbar{\gamma}|_{\gamma>0}} =  \frac{-\zeta\ + \frac{\sigma_{\alpha}}{4} \frac{\phi(\eta)}{1-\Phi(\eta)}}{ \zeta\ +  \frac{\sigma_{\alpha}}{4}  \frac{\phi(\eta)}{\Phi(\eta)}} \equiv \upsilon.
\end{equation}
The average duration of rare transitions is thus given by $T_e = (1 + \overbar{\Lambda}_+/\overbar{\Lambda}_- ) T_{\gamma<0}$. Dividing over the total length of time $T$, gives the final result for the probability of rare transitions:
\begin{equation}
\prob_r  \simeq  (1 + \overbar{\Lambda}_+/\overbar{\Lambda}_- ) \prob(\gamma<0) =  (1 + \upsilon)(1 - \Phi(\eta) ).
\end{equation}

\subsection{Summary of the analytical results}

With the analysis above, we synthesize the results by a total probability law argument:\begin{multline}\label{eq:mathanalyic}
\pdf_x( x ) = (1-\prob_r) \sqrt{\frac{\overbar \gamma |_{\gamma>0}}{\pi \sigma_F^2 }} \exp \biggl( \frac{-\overbar \gamma |_{\gamma>0}}{\sigma_F^2} x^2 \biggr) \\+ 
\prob_r  \frac{ \sqrt{ 2 \pi} \overbar\gamma|_{\gamma>0}
}{\sigma_F^2 \sigma_{\alpha}\omega_0\ \overbar  T^2 (1-\Phi(\eta))   }  \int\limits_\zeta^{\abs{x}}\!\int\limits_0^\infty
\frac{\log(\abs{x}/\xi_0)}{y^2(\abs{x}/(\xi_0-\zeta))} \\ \times\exp\biggl(-8\frac{(y
+ \zeta \omega_0)^2}{ \sigma_{\alpha}^2\omega_0^2}  - \frac{\pi}{4 \overbar
T^2 y^2 } \log(\abs{x}/\xi_0)^2 - \frac{\overbar\gamma|_{\gamma>0}}{\sigma_F^2}
(\xi_0-\zeta)^2 \bigg) \,dy\,d\xi_0,
\end{multline}
where the integral is   zero for $\abs{x} < \zeta$. Note that $\sigma_F^2 = \pi S_F(\omega_0)/\omega_0^2$  for the case of a broadbanded additive   excitation and $\sigma_F^2 = \nu^2/2\omega_0^2$ for the case of a white noise additive excitation   with intensity $\nu$ (i.e. $F(t) = \nu \dot W(t)$). 

\begin{figure}[htb]
    \centering
                \includegraphics[width=\textwidth]{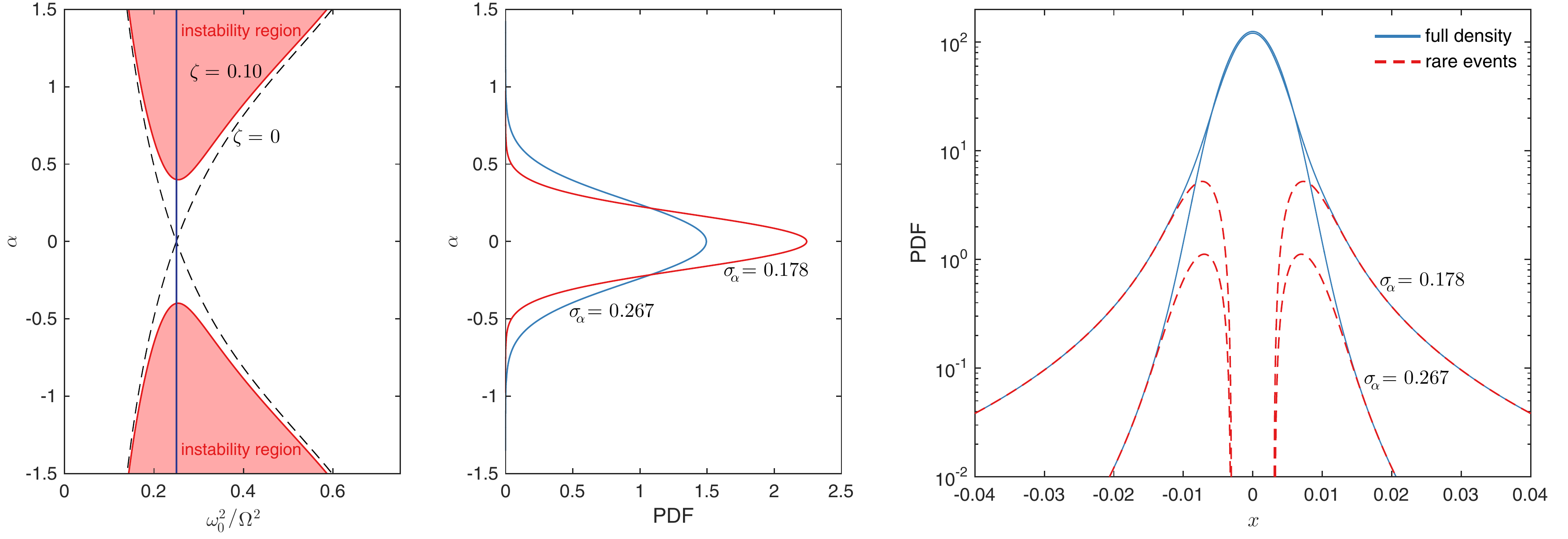}
    \caption{(Left) Stability diagram for constant (in time) $\alpha$ and for finite damping $\zeta=0.1$ (red shaded region); (Middle) pdf of the excitation $\alpha(t)$ for two different values of $\sigma_{\alpha}$; (Right) The corresponding analytical pdf of Mathieu's equation according to the decomposition-sythesis method~\labelcref{eq:mathanalyic}, where we have  highlighted the conditionally rare event component  (dashed red) in the full pdf (solid blue),  showing that the tails are completely described by rare events, whereas the background dynamics fully determine the core  (the response pdf is shown assuming $\ell_\alpha = 10.0$,   other parameters are described in~\cref{sec:comps}).}
    \label{fig:superpos}
\end{figure}

We emphasize that~\cref{eq:mathanalyic} is a heavy-tailed and symmetric  probability measure (non-negative and integrates to one). In~\cref{fig:superpos} we  present the stability diagram for constant (in time) $\alpha$ and for finite damping $\zeta=0.1$ (red shaded region in the left plot). For the case where  $\Omega=2\omega_0$ we apply a random amplitude parametric excitation with pdf for $\alpha(t)$ shown in the middle plot for two different values of $\sigma_{\alpha}$. (Note that there are other important factors that play an important role in the form of the tails, not captured in the stability diagram or the pdf of $\alpha(t)$, such as the correlation function of the process $\alpha(t)$.) In the right plot,  the corresponding response pdf as computed through~\cref{eq:mathanalyic} for the same two values of $\sigma_{\alpha}$ are shown, illustrating the heavy-tailed component, due to the conditionally rare dynamics, and the core of the pdf, due to the conditionally background state. It is clear that   transient instabilities, rare responses,     fully determine the tails of the pdf, while the background dynamics  specify    the core of the   distribution, but contribute essentially zero probability to the tails.

\section{Comparisons with Monte-Carlo simulations}\label{sec:comps}

To illustrate the accuracy of our approximation~\labelcref{eq:mathanalyic},  we compare the analytical formula obtained via the probabilistic-decomposition technique with direct Monte-Carlo simulations of the original system. To perform comparisons we use a unit  white noise $\dot W(t)$ process with intensity $\nu$ for the additive forcing term $F(t)$, and furthermore non-dimensionalize time by $\omega_0$ in~\cref{eq:problemstate_stoch_Mathieu} so that, 
\begin{equation}\label{eq:referenceEqMC}
\ddot x(t) + 2   \zeta   \dot x(t) +  ( 1 +   \alpha(t)  \sin  2   t )  x(t) =     \nu  \dot W(t).
\end{equation}
Thus, for this additive forcing we have $\sigma_F^2 = \nu^2/2\omega_0^2$.

To perform the Monte-Carlo simulations, we compute  $3000$ realization of the above equation  using the Euler-Maruyama method  with  time step $dt =  5 \times 10^{-3}$ from $t= 0$ to $t = 5500$ and discard the first $500$ time units, to ensure a statistical steady state. Moreover, we generate realizations of  the stochastic process $\alpha(t)$ directly from the autocorrelation function by a statistically exact  and efficient method using the  circulant embedding technique~\cite{2013arXiv1308.0399K}.

\begin{figure}[htb]
    \centering
    \includegraphics[width=\textwidth]{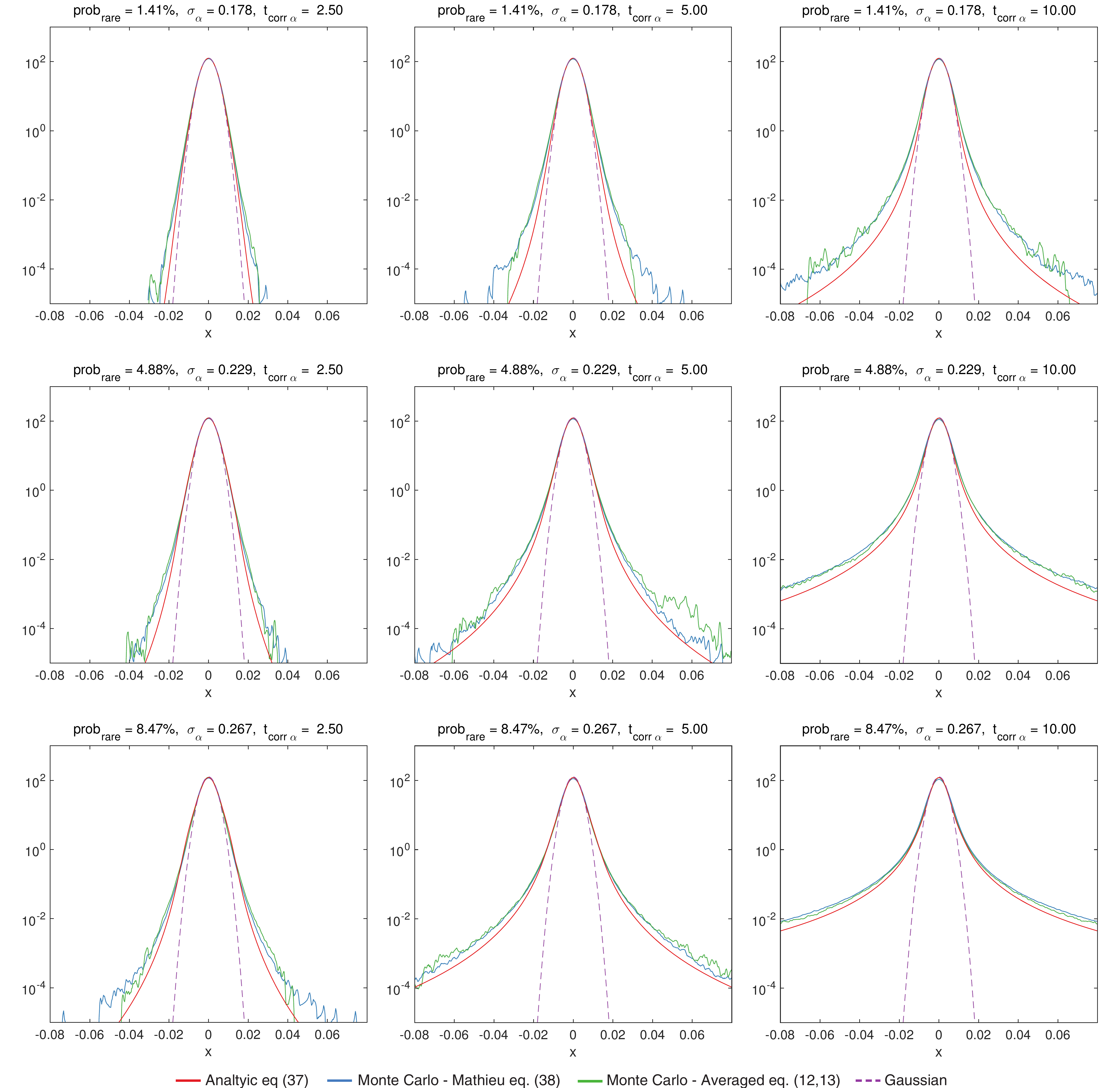}
    \caption{Comparisons of Monte-Carlo simulations (blue curve) of Mathieu's~\cref{eq:referenceEqMC}  the `true' density, the averaged system~\cref{eq:averagedeqnsx1,eq:averagedeqnsx2} (green curve), and  our analytical formula~\labelcref{eq:mathanalyic} (red curve), along with a  Gaussian approximation (dashed purple curve), where the variance is obtained by a traditional averaging procedure. The left column corresponds to a parametric excitation process with correlation time-scale of $\ell_\alpha = 2.5$, the middle column $\ell_\alpha = 5.0$, and the right column $\ell_\alpha = 10.0$; note the time scale of damping is $t_\text{damp}\sim 1/\zeta = 10.0$. In the top row the excitation standard deviation is $\sigma_\alpha = 0.148$, middle row $\sigma_\alpha = 0.229$, and the bottom row $\sigma_\alpha = 0.267$. Fixed system parameters: $\zeta = 0.01$ and $\nu = 0.002$.}
    \label{fig:results}
\end{figure}

In~\cref{fig:results}   we show nine cases of varying intermittency levels. We fix the system parameters $\zeta  = 0.1$ and  $\nu\ = 0.002$. In the figure, we show the response pdf for  three different correlation times $\ell_\alpha = 2.5,\,5.0,\,10.0$ of    and for three different values of the standard deviation of the parametric excitation $\alpha(t)$: $\sigma_\alpha=0.178,\,0.229,\,0.267$; For the least  intermittent case    the standard deviation of $\alpha(t)$ is $\sigma_\alpha = 0.178$ so that  rare event transitions occur with probability $\prob_r = 0.0141$, for  $\sigma_\alpha = 0.229$ with probability $\prob_r = 0.0488$, and for the most intermittent case $\prob_r = 0.0847$. 

Overall, the results presented, along with  additional numerical comparisons, show good  quantitative agreement for both the tails and the core of the distribution  between our analytic formula in~\cref{eq:mathanalyic} and the `true' density from Monte-Carlo simulations; indeed, the quantitative agreement between our formula and the actual density  is  found to be robust  across a range of parameters that satisfy the assumptions. We observe our result performs better with the averaged system, on which we directly derived the density, than the original system for cases with larger correlation times $\ell_\alpha$ and larger variance $\sigma_\alpha^2$ in the parametric excitation process $\alpha(t)$; this behavior is expected  since averaging introduces well known errors that increase for  more intermittent regimes, because   the instabilities in such cases  lead to even larger amplitude responses. Moreover, we observe that even in extremely intermittent regimes, where our  assumptions start to be violated, namely the statistical independence of rare events, our analytical formula is  still able to  capture the asymptotic behavior of the tails.

Clearly, the results show that the response pdf is far from Gaussian, and highlight the    non-Gaussianity of roll motions during parametric resonance. Despite the simplicity and theoretical nature of our model, the overall characteristics and shape of the distribution is the same as that observed in Monte-Carlo simulations using advanced  hydrodynamic codes with three degrees of freedom  that account for coupled pitch, heave, and roll motions~\cite{vadim2011head,belenkyprobprop}.

\section{Conclusions}\label{sec:conclusions}

In this work we derived an analytical approximation to the heavy-tailed stationary measure of Mathieu's equation under  parametric excitation by a correlated stochastic process  in a regime undergoing  intermittent parametric resonances.   We derived the formula for the case when the spectrum of the noise is peaked at the main resonant frequency. To derive the pdf for the response we averaged the governing equation over the fast frequency to arrive at a set of parametrically excited   processes that govern the slow dynamics. We then  applied the  probabilistic decomposition-synthesis method to the slow variables.  We demonstrated the accuracy of the final   formula for the pdf  through direct comparisons with Monte-Carlo results for a range parameters that influence the rare event transition probability level and severity of the resonance phenomena; the analytical formula showed excellent quantitative agreement with results from numerical simulations across a wide range of intermittency levels. The approach is also directly applicable for the determination of the local maxima of the response.

The presented analysis   paves the way for the analytical treatment of more realistic ship roll models.  Future work  includes the inclusion of nonlinear terms, in particular, softening nonlinearity in the restoring force, which would modify  the   statistical characteristics of the tails. Such a problem could   be considered through the current framework by analytical modeling of the nonlinear terms during the rare transitions, in combination with appropriate nonlinear closures for the background stochastic attractor. Additional work will include application of the decomposition-synthesis method in a data-driven context using models containing model error to derive tail estimates.

\section*{Acknowledgments}
This research has been partially supported  by the Office of Naval Research (ONR) Grant ONR N00014-14-1-0520 and the Naval Engineering Education
Center (NEEC) Grant 3002883706. We thank  Dr. Vadim Belenky for numerous stimulating discussions. We are also  grateful to Profs. Francescutto, Neves, and Vassalos for the invitation  to prepare a manuscript for this special issue of Ocean Engineering Journal on Stability and Safety of Ships and Ocean Vehicles.

\clearpage
\bibliography{main,library}

\end{document}